\documentclass[ECP]{ejpecp} 

\SHORTTITLE{Chung's LIL for SPDEs}

\TITLE{Chung's Law of the Iterated Logarithm for a Class of Stochastic Heat Equations}

\AUTHORS{
  Jiaming Chen\footnote{Dept. of Mathematics, University of Rochester, Rochester, NY 14627,
    \EMAIL{jchen143@ur.rochester.edu}}}
  

\AMSSUBJ{60H15} 
\AMSSUBJSECONDARY{60G17} 

\SUBMITTED{February 28, 2023} 
\ACCEPTED{August 27, 2023} 

\VOLUME{0}
\YEAR{2020}
\PAPERNUM{0}
\DOI{10.1214/YY-TN}

\ABSTRACT{We establish a Chung-type law of the iterated logarithm for the solutions of a class of stochastic heat equations driven by a multiplicative noise whose coefficient depends on the solution, and this dependence takes us away from Gaussian setting. Based on the literature on small ball probabilities and the technique of freezing coefficients, the limiting constant in Chung’s law of the iterated logarithm can be evaluated almost surely.}

\newcommand{\R}{\mathbb{R}}

\begin{document}



\section{Introduction}
Let $\{B_t\}_{t\geq 0}$ be a standard one-dimensional Brownian motion. Chung \cite{chung1948maximum} proved that, with probability one,
\[
\liminf_{t\to\infty}\left(\frac{\log\log t}{t}\right)^{\frac{1}{2}}\sup_{0\leq s\leq t}\vert B_s\vert =\frac{\pi}{\sqrt{8}}.
\]

Generally speaking, Chung's law of the iterated logarithm (Chung's LIL) characterizes the lower envelope $(\liminf)$ for the local oscillations of the sample path.

Chung's LIL has already been applied to empirical processes \cite{mogul1980law}, some stochastic integrals \cite{remillard1994chung}, the functional form of Chung's LIL for Brownian motion \cite{kuelbs1994lim}, Brownian sheet \cite{talagrand1994small}, the iterated Brownian motion \cite{khoshnevisan1996chung} and a hypoelliptic Brownian motion in Heisenberg group \cite{carfagnini2022small}. Furthermore, Luan and Xiao \cite{luan2010chung} establish Chung's LIL for a class of anisotropic Gaussian random fields with stationary increments, which is applicable to space-time Gaussian random fields and a solution of stochastic partial differential equations (SPDEs). Recently, Lee and Xiao \cite{lee2023chung} proved Chung's LIL in a wider class of Gaussian random fields that may not necessarily have stationary increments. However, we still limited ourselves to Gaussian processes in most of the results mentioned above, and that there are not many explorations of non-Gaussian processes.

The purpose of this paper is to study Chung's LIL for a class of non-Gaussian processes at the origin. This is mainly motivated by recent papers \cite{lee2023chung} and \cite{athreya2021small} since it is well-known that a key step in establishing Chung's LIL is small ball probability. Small ball problems have been well studied and have a long history, and one can see \cite{li2001gaussian} for an overview of known results on Gaussian processes and references on other processes. In short, we investigate the behavior of the probability that a stochastic process $X_t$ starting at the origin stays near its initial location for a long time, i.e.,
\[P\left(\sup_{0\leq t\leq T}\vert X_t\vert<\varepsilon\right)\]
where $\varepsilon>0$ is small. When small ball probabilities are available, it is not too difficult to derive Chung’s LIL in Gaussian cases, and one major step is to apply the second Borel-Cantelli lemma to a sequence of independent events. However, particularly in such a non-Gaussian process, obtaining Chung’s LIL is a completely different matter.

\section{Main Result}
We consider a class of stochastic heat equations given by 
\begin{equation}
\label{heat}
\begin{split}
\partial_t &u(t,x)=\frac{1}{2}\partial_x^2 u(t,x)+\sigma(u(t,x))\dot{W}(t,x),\\
&u(0,x)=u_0(x)\equiv 0,
\end{split}
\end{equation}
on the circle with $x\in[0,1]$ and endpoints identified, $t\in\R_+$, and $\dot{W}(t,x)$ is space-time white noise. We assume that the function $\sigma:\R\to \R$ is globally Lipschitz continuous.
\begin{hypothesis}
\label{hypothesis}
There exists a constant $\mathcal{D}\geq 0$ such that, for all $u,v\in\R$,
\begin{equation*}
\vert \sigma(u)-\sigma(v)\vert\leq\mathcal{D}|u-v|.
\end{equation*}
\end{hypothesis}

In fact, \eqref{heat} is not well-posed since the solution $u$ is not differentiable and $\dot{W}$ exists as a generalized function. However, under the Hypothesis \ref{hypothesis} and $u_0(x)\equiv 0$, we can define the mild solution $u(t,x)$ to \eqref{heat} in the sense of Walsh \cite{walsh1986anintroductiontostochastic}:
\begin{equation}
\label{sol}
u(t,x)=\int_0^1 G(t,x-y)u_0(y)dy+\int_0^t\int_0^1G(t-s,x-y)\sigma(u(s,y))W(dyds),
\end{equation}
where $G:\R_+\times[0,1]\to\R$ is the fundamental solution of the heat equation
\begin{equation*}
\begin{split}
\partial_t&G(t,x)=\frac{1}{2}\partial^2_xG(t,x),\\
&G(0,x)=\delta_0(x).
\end{split}
\end{equation*}

We are now ready to state the main result of this paper.
\begin{theorem}
\label{thm}
Under the Hypothesis \ref{hypothesis}, if $u(t,x)$ is the solution to \eqref{heat} with $u_0(x)\equiv 0$ and $\sigma(0)\neq 0$, then there exists a non-degenerate constant $\kappa\in(0,\infty)$ such that we have, almost surely,
\begin{equation}
\label{thmeq1}
\liminf_{r\to 0^+}\sup_{\substack{0\leq t\leq r^4\\ x\in\left[0,r^2\right]}}\frac{\vert u(t,x)\vert}{f(r)}=\kappa,
\end{equation}
where $f(r):=r\left[\log\log\left(\frac{1}{r}\right)\right]^{-\frac{1}{6}}$. 
\end{theorem}

\begin{remark}
When $\mathcal{D}=0$, \eqref{heat} becomes a linear stochastic heat equation with a constant coefficient, and the result is indicated in \cite[Theorem 7.4]{lee2023chung}.
\end{remark}

Let us summarize major steps in this paper. We truncate the coefficient $\sigma$ and determine the probability that the solution remains unchanged despite the truncation. This probability can then be bounded above by a tail bound in \cite{athreya2021small}. With the first Borel-Cantelli lemma, we demonstrate that this truncation does not affect the limiting term almost surely. 

The solution dependent coefficient displaces us from the Gaussian setting, rendering it unfeasible to identify a sequence of independent events due to this dependency. Nevertheless, we could freeze the truncated coefficient and construct a sequence of events involving a Gaussian process, and thus we could use the estimation in \cite{lee2023chung} and approximate the error terms using the freezing technique in \cite{athreya2021small}. To be precise, given the result of the limiting term for a linear SPDE with a constant coefficient in \cite{lee2023chung}, we may approximate for a non-Gaussian process if we can properly estimate the limiting term with respect to their differences. This idea first appeared in \cite{pospivsil2007parameter} for the stochastic heat equation, yet it was used to perform $L^2$ estimates instead of almost sure estimates. Lastly, we apply the first Borel-Cantelli lemma again along with this freezing technique to derive the limiting constant.

The rest of the article is organized as follows. In Section \ref{preliminary}, we provide some preliminaries in the literature to prove Theorem \ref{thm}. In Section \ref{proofthm}, we conclude the paper by showing that Theorem \ref{thm} can be derived from a Gaussian result.

\section{Preliminary}
\label{preliminary}
In the paper \cite{athreya2021small}, the authors discovered a crucial lemma and a powerful remark that give critical tail bounds on the noise term.

\begin{lemma}\cite[Lemma 3.4]{athreya2021small} \label{ajmlemma}
Assume that there is a positive constant $\mathcal{C}$ such that $\sigma(u)<\mathcal{C}$ for all $u\in\R$, then there exist positive constants $\textbf{K}_1$, $\textbf{K}_2$ such that, for all $\varepsilon,\lambda>0$, we have
\begin{equation}
\label{ajmlemmaeq}
P\left(\sup_{\substack{0\leq t\leq \varepsilon^4\\ x\in\left[0,\varepsilon^2\right]}}\vert u(t,x)\vert>\lambda \varepsilon\right)\leq \textbf{K}_1\exp\left(-\frac{\textbf{K}_2\lambda^2}{\mathcal{C}^2}\right).
\end{equation}
\end{lemma}

\begin{remark}
\label{ajmrem}
One can obtain \eqref{ajmlemmaeq} by letting $\alpha=1$ and $u_0\equiv 0$ in \cite[Lemma 3.4]{athreya2021small}. It was also pointed out in \cite[Remark 3.1]{athreya2021small} that if $\vert\sigma(u(t,x))\vert\leq C_1\varepsilon$, then one can replace the right-hand side of \eqref{ajmlemmaeq} with $\textbf{K}_1\exp\left(-\textbf{K}_2\frac{\lambda^2}{C_1^2\varepsilon^2}\right).$
\end{remark}

Recent paper \cite{lee2023chung} established a general framework that is useful for studying the regularity properties of sample functions of anisotropic Gaussian random fields and can be directly applied to the solution of linear SPDEs.

The following Theorem provides us with a preliminary result in a Gaussian process, and we can utilize a freezing technique to connect this result to a non-Gaussian process.

\begin{theorem}\cite[Theorem 7.4]{lee2023chung}\label{lxthm} Suppose that $u(t,x)$ is the solution to \eqref{heat} with $\sigma\equiv 1$, then there is a constant $\kappa_0\in(0,\infty)$ such that we have, almost surely,
\begin{equation*}
\liminf_{r\to 0^+}\sup_{\substack{0\leq t\leq r^4\\ x\in\left[0,r^2\right]}}\frac{\vert u(t,x)\vert}{f(r)}=\kappa_0,
\end{equation*}
where $f(r)=r\left[\log\log\left(\frac{1}{r}\right)\right]^{-\frac{1}{6}}$. Moreover, $\kappa_0$ coincides with a limit in \cite[Theorem 4.4]{lee2023chung}, which is called the small ball constant of $u$ on the region $0\leq t\leq r^4,x\in\left[0,r^2\right]$, if it exists.
\end{theorem}

\begin{remark}
Here we estimate the limiting term in a rectangle region, which is compatible with the metric in \cite{lee2023chung}. Furthermore,  \cite{dalang2017polarity} have established a harmonizable representation for the solution of linear stochastic heat equations with constant coefficients, which implies that $Q=6$ in \cite[Theorem 7.4]{lee2023chung} when $\dot{W}$ is space-time white noise.
\end{remark}

\section{Proof of Theorem \ref{thm}}
\label{proofthm}
Given that $\sigma(0)\neq 0$, we let positive constant $m=2\vert\sigma(0)\vert$ and define a truncated function
\begin{equation}
\label{sigmatilde}
\tilde{\sigma}(x)=\begin{cases}
\sigma(x),&\text{if}~\vert\sigma(x)\vert\leq m;\\
m,&\text{if}~\sigma(x)> m;\\
-m,&\text{if}~\sigma(x)<-m.\\
\end{cases}
\end{equation}
Clearly, $\tilde{\sigma}$ is a bounded function and we now start showing that the limiting term in \eqref{thmeq1} can be approximated by
\begin{equation*}
\liminf_{r\to 0^+}\sup_{\substack{0\leq t\leq r^4\\ x\in\left[0,r^2\right]}}\frac{\vert \tilde{u}(t,x)\vert}{f(r)},
\end{equation*}
where $\tilde{u}$ solves the stochastic heat equation
\begin{equation}
\label{tildeudef}
\begin{split}
\partial_t &\tilde{u}(t,x)=\frac{1}{2}\partial_x^2 \tilde{u}(t,x)+\tilde{\sigma}(\tilde{u}(t,x))\dot{W}(t,x),\\
&\tilde{u}(0,x)\equiv 0.
\end{split}
\end{equation}

Let $r_n=a^{-n}$ with $a>1$, and then from the definition of $\tilde{\sigma}$ in \eqref{sigmatilde} and the uniqueness property of the solution of the stochastic heat equation, we have
\begin{equation}
\label{tildeu}
P\left(u(t,x)\neq \tilde{u}(t,x)~\text{for some}~(t,x)\in\left[0,r_n^4\right]\times\left[0,r_n^2\right]\right)\leq P\left(\sup_{\substack{0\leq t\leq r_n^4\\ x\in\left[0,r_n^2\right]}}\vert \sigma(\tilde{u}(t,x))\vert>m\right).
\end{equation}

If $\mathcal{D}=0$ $[$i.e., $\sigma$ is constant$]$, then the above probability in \eqref{tildeu} is zero because $\vert\sigma(\tilde{u}(t,x))\vert=\vert\sigma(0)\vert=\frac{m}{2}<m$ for all $(t,x)\in\left[0,r_n^4\right]\times\left[0,r_n^2\right]$. On the other hand, we assume that $\mathcal{D}>0$ and Hypothesis \ref{hypothesis} implies,
\begin{equation*}
\mathcal{D}\vert\tilde{u}(t,x)\vert\geq \vert\sigma(\tilde{u}(t,x))-\sigma(0)\vert\geq \vert\sigma(\tilde{u}(t,x))\vert-\frac{m}{2}.
\end{equation*}
We recall that $\vert\tilde{\sigma}(x)\vert\leq m$ in \eqref{sigmatilde}, so \eqref{tildeu} and Remark \ref{ajmrem} together yield
\begin{equation*}
\begin{split}
P\left(u(t,x)\neq \tilde{u}(t,x)~\text{for some}~(t,x)\in\left[0,r_n^4\right]\times\left[0,r_n^2\right]\right)&\leq P\left(\sup_{\substack{0\leq t\leq r_n^4\\ x\in\left[0,r_n^2\right]}}\vert \tilde{u}(t,x)\vert>\frac{m}{2\mathcal{D}}\right)\\
&\leq \textbf{K}_1\exp\left(-\frac{\textbf{K}_2m^2}{4\mathcal{D}^2m^2r_n^2}\right)\\
&=\textbf{K}_1\exp\left(-\frac{\textbf{K}_2}{4\mathcal{D}^2r_n^2}\right).
\end{split}
\end{equation*}
The first Borel-Cantelli lemma implies that, with probability one,
\begin{equation*}
\lim_{r\to 0^+}\sup_{\substack{0\leq t\leq r^4\\ x\in\left[0,r^2\right]}}\frac{\vert u(t,x)-\tilde{u}(t,x)\vert}{f(r)}=0.
\end{equation*}
Hence, we conclude that, with probability one,
\begin{equation}
\label{uutilde}
\liminf_{r\to 0^+}\sup_{\substack{0\leq t\leq r^4\\ x\in\left[0,r^2\right]}}\frac{\vert u(t,x)\vert}{f(r)}=\liminf_{r\to 0^+}\sup_{\substack{0\leq t\leq r^4\\ x\in\left[0,r^2\right]}}\frac{\vert \tilde{u}(t,x)\vert}{f(r)}.
\end{equation}

Now it suffices to prove Theorem \ref{thm} for $\tilde{u}$ that solves the stochastic heat equation \eqref{tildeudef} driven by a multiplicative noise whose coefficient is bounded.

As we previously mentioned in the introduction, it may not be feasible to directly show Theorem \ref{thm} using the Borel-Cantelli lemma because of the solution dependent coefficient. Instead of seeking for a sequence of independent events, we freeze the coefficient $\tilde{\sigma}(\tilde{u})$ with $\tilde{\sigma}(\tilde{u}_0)$, that is, we may approximate the solution $\tilde{u}(t,x)$ by a Gaussian random variable, at least in a random time region. 

We let $\tilde{u}_g(t,x)$ satisfy the stochastic heat equation
\begin{equation*}
\begin{split}
\partial_t &\tilde{u}_g(t,x)=\frac{1}{2}\partial_x^2 \tilde{u}_g(t,x)+\tilde{\sigma}(\tilde{u}(0,x))\dot{W}(t,x),\\
&\tilde{u}_g(0,x)\equiv 0,
\end{split}
\end{equation*}
and then $\tilde{u}_g(t,x)$ is the solution of a linear SPDE with a coefficient $\tilde{\sigma}(0)=\sigma(0)$. The following statement evaluates the sum of the tail probability of the difference between $\tilde{u}(t,x)$ and $\tilde{u}_g(t,x)$, and this difference can be well controlled if the time region is suitably chosen. 

From the definition of the mild solution \eqref{sol}, we denote
\begin{equation*}
\begin{split}
D(t,x)&=\tilde{u}(t,x)-\tilde{u}_g(t,x)\\
&=\int_{[0,1]\times[0,t]}G(t-s,x-y)[\tilde{\sigma}(\tilde{u}(s,y))-\tilde{\sigma}(\tilde{u}(0,y))]W(dyds),
\end{split}
\end{equation*}
and 
\begin{equation*}
\widetilde{D}_n(t,x)=\int_{[0,1]\times[0,t]}G(t-s,x-y)\left[\tilde{\sigma}\left(\tilde{u}\left(s\wedge\tau_n^g,y\right)\right)-\tilde{\sigma}(\tilde{u}(0,y))\right]W(dyds),
\end{equation*}
where $\tau_n^g$ is a stopping time defined as
\begin{equation*}
\tau^g_n=\inf\left\lbrace t>0:\vert \tilde{u}(t,x)-\tilde{u}(0,x)\vert>g(r_{n})\text{~for some $x\in\left[0,r_n^2\right]$}\right\rbrace.
\end{equation*}
Clearly, on the event $F_n=\left\lbrace\tau_n^g\geq r_n^4\right\rbrace$, we have
\begin{equation*}
\sup_{\substack{0\leq t\leq r_n^4\\ x\in\left[0,r_n^2\right]}}\vert D(t,x)\vert=\sup_{\substack{0\leq t\leq r_n^4\\ x\in\left[0,r_n^2\right]}}\left\vert \widetilde{D}_n(t,x)\right\vert,
\end{equation*}
which leads to, for every $0<\varepsilon<1$,
\begin{equation}
\label{ineq}
P\left(\sup_{\substack{0\leq t\leq r_n^4\\ x\in\left[0,r_n^2\right]}}\vert D(t,x)\vert>r_n^{1+\varepsilon}\right)\leq P\left(\sup_{\substack{0\leq t\leq r_n^4\\ x\in\left[0,r_n^2\right]}}\left\vert \widetilde{D}_n(t,x)\right\vert>r_n^{1+\varepsilon}\right)+P(F_n^c).
\end{equation} 
Choosing $g(r_n)=r_n^{\frac{1+\varepsilon}{2}}$, we can apply Hypothesis \ref{hypothesis}, Lemma \ref{ajmlemma} and Remark \ref{ajmrem} to obtain
\begin{equation}
\label{tildeD}
\begin{split}
P\left(\sup_{\substack{0\leq t\leq r_n^4\\ x\in\left[0,r_n^2\right]}}\left\vert \widetilde{D}_n(t,x)\right\vert>r_n^{1+\varepsilon}\right)&\leq \textbf{K}_1\exp\left(-\frac{\textbf{K}_2r_n^{2\varepsilon}}{\mathcal{D}^2r_n^{1+\varepsilon}}\right)\\
&=\textbf{K}_1\exp\left(-\frac{\textbf{K}_2}{\mathcal{D}^2r_n^{1-\varepsilon}}\right),
\end{split}
\end{equation}
and
\begin{equation}
\label{Fc}
\begin{split}
P(F_n^c)&\leq P\left(\sup_{\substack{0\leq t\leq r_n^4\\ x\in\left[0,r_n^2\right]}}\vert \tilde{u}(t,x)-\tilde{u}(0,x)\vert>r_n^{\frac{1+\varepsilon}{2}}\right)\\
&\leq \textbf{K}_1\exp\left(-\frac{\textbf{K}_2r_n^{1+\varepsilon}}{r_n^2m^2}\right)\\
&=\textbf{K}_1\exp\left(-\frac{\textbf{K}_2}{r_n^{1-\varepsilon}m^2}\right).
\end{split}
\end{equation}

We recall that $r_n=a^{-n}$, $a>1$ and $0<\varepsilon<1$, so summing up probabilities in \eqref{tildeD} and \eqref{Fc} yield
\begin{equation*}
\sum_{n=1}^\infty P\left(\sup_{\substack{0\leq t\leq r_n^4\\ x\in\left[0,r_n^2\right]}}\left\vert \widetilde{D}_n(t,x)\right\vert>r_n^{1+\varepsilon}\right)<\infty,
\end{equation*}
and 
\begin{equation*}
\sum_{n=1}^\infty P(F^c_n)<\infty.
\end{equation*}
According to \eqref{ineq}, the sum of the following sequence of events
\begin{equation*}
\sum_{n=1}^\infty P\left(\sup_{\substack{0\leq t\leq r_n^4\\ x\in\left[0,r_n^2\right]}}\frac{\vert D(t,x)\vert}{f(r_n)}>r_n^{\varepsilon}\left[\log\log\left(\frac{1}{r_n}\right)\right]^{\frac{1}{6}}\right)=\sum_{n=1}^\infty P\left(\sup_{\substack{0\leq t\leq r_n^4\\ x\in\left[0,r_n^2\right]}}\vert D(t,x)\vert>r_n^{1+\varepsilon}\right)<\infty,
\end{equation*}
with the first Borel-Cantelli lemma shows that, almost surely,
\begin{equation*}
\lim_{r\to0^+}\sup_{\substack{0\leq t\leq r^4\\ x\in\left[0,r^2\right]}}\frac{\vert \tilde{u}(t,x)-\tilde{u}_g(t,x)\vert}{f(r)}=0,
\end{equation*}
and
\begin{equation}
\label{uug}
\liminf_{r\to 0^+}\sup_{\substack{0\leq t\leq r^4\\ x\in\left[0,r^2\right]}}\frac{\vert \tilde{u}(t,x)\vert}{f(r)}=\liminf_{r\to 0^+}\sup_{\substack{0\leq t\leq r^4\\ x\in\left[0,r^2\right]}}\frac{\vert \tilde{u}_g(t,x)\vert}{f(r)}.
\end{equation}

We can write $\tilde{u}_g(t,x)=\sigma(0)v(t,x)$, where $v(t,x)$ solves the stochastic heat equation
\begin{equation*}
\begin{split}
\partial_t &v(t,x)=\frac{1}{2}\partial_x^2 v(t,x)+\dot{W}(t,x),\\
&v(0,x)\equiv 0,
\end{split}
\end{equation*}
and immediately get
\begin{equation}
\label{v}
\liminf_{r\to 0^+}\sup_{\substack{0\leq t\leq r^4\\ x\in\left[0,r^2\right]}}\frac{\vert \tilde{u}_g(t,x)\vert}{f(r)}=\sigma(0)\cdot\liminf_{r\to 0^+}\sup_{\substack{0\leq t\leq r^4\\ x\in\left[0,r^2\right]}}\frac{\vert v(t,x)\vert}{f(r)}.
\end{equation}
By Theorem \ref{lxthm}, \eqref{uutilde}, \eqref{uug} and \eqref{v}, we end up with, almost surely,
\[
\liminf_{r\to 0^+}\sup_{\substack{0\leq t\leq r^4\\ x\in\left[0,r^2\right]}}\frac{\vert u(t,x)\vert}{f(r)}=\sigma(0)\kappa_0,
\]
which completes the proof with $\kappa=\sigma(0)\kappa_0$ in Theorem \ref{thm}. 




\begin{acks}
The author would like to thank his advisor Professor Carl Mueller and the anonymous referees for constructive comments.
\end{acks}


\end{document}